\documentclass[11pt]{amsart}

\usepackage{amsfonts}
\usepackage{amssymb}
\usepackage{mathrsfs}
\usepackage{amsrefs}

\theoremstyle{plain}
\newtheorem{thm}{Theorem}[section]
\newtheorem{prop}[thm]{Proposition}

\newtheorem{lemma}[thm]{Lemma}

\renewcommand{\mod}[1]{\left\lvert #1 \right\rvert}

\newcommand{\norm}[1]{\left\| #1 \right\|}
\newcommand{\cl}[1]{\overline{#1}}
\newcommand{\bb}[1]{\mathbb{#1}}
\newcommand{\cal}[1]{\mathcal{#1}}

\DeclareMathOperator{\wk}{{weak}^\ast}

\begin{document}
\title[Real Hardy spaces]{Function Theory in Real Hardy Spaces}
\author[M. Raghupathi]{Mrinal Raghupathi}
\address{Department of Mathematics, University of Houston, Houston, TX, 77204}
\thanks{The first author was supported by NSF grant DMS 0300128.}
\email{mrinal@math.uh.edu}
\urladdr{http://www.math.uh.edu/~mrinal}

\author[D. Singh]{Dinesh Singh}
\address{Department of Mathematics, University of Delhi, Delhi, India}
\email{dineshsingh1@gmail.com}
\thanks{The second author would like to acknowledge the Department of Mathematics, University of Houston and the Mathematical Sciences Foundation, Delhi}
\begin{abstract}
 We show that many classical results in Hardy space theory have exact analogues when the Fourier coefficients are allowed only to be real. 
\end{abstract}

\maketitle

\section{Introduction}\label{intro}

In \cite{JR} and \cite{MS} the non-trivial closed shift invariant subspaces of the real linear space of $L^1$ functions with real Fourier coefficients are classified. Their structure is analogous to the classical setting. We will show here that a much simpler proof is possible which exploits the classical situation. Our technique will allow us not only to give short proofs of the results in \cite{MS} and \cite{JR} but will lead to the analogue of many other classical function theorietic results. In the recent literature there has been some interest in the real disk algebra. We point of the work of Wick \cite{wick}, which deals with the notion of Bass stable rank for the real disk algebra.

In section \ref{prelim} we make some elementary definitions and lay down notation. We will describe the complexification of a real vector space and a projection operator. After deriving some basic properties of this projection we will proceed, in the remaining sections, to derive the function theoretic results in this new setting. In section \ref{misc} we will describe how these ideas can be extended to the case of compact abelian groups with ordered duals and to functions on the line with real Fourier transforms, where we can prove the analogue of Wiener's and Lax's theorem.

For the classical results we refer the reader to the books of Helson \cite{HH} and Hoffman \cite{KH}.

\section{Preliminaries}\label{prelim}

Let $G$ be a locally compact abelian group and let $m$ denote Haar measure on the group. In the case where $G$ is compact we will assume that $m$ has been normalized so that $m(G)=1$. We will denote by $\hat{G}$ the set of continuous homomorphisms of $G$ into the circle group $\bb{T}$. The dual group is given the topology of uniform convergence on compact sets. i.e. a net $\chi_t \to \chi$ if and only if for every compact subset $K\subseteq G$, $\chi_t|_K \to \chi|_K$ uniformly on $K$. It is a fact that if $G$ is compact, then $\hat{G}$ is discrete and if $G$ is discrete, then $\hat{G}$ is compact. It is well known that $G$ naturally embeds in $\hat{\hat{G}}$ via the map $g\mapsto \delta_g$, where $\delta_g(\chi) = \chi(g)$. The Pontryagin duality theorem asserts that this embedding is in fact an isomorphism of topological groups. It follows that $G$ is compact if and only if $\hat{G}$ is discrete and that $G$ is discrete if and only if $\hat{G}$ is compact. 

For $1\leq p\leq \infty$, let $L^p(G)$ denote the Lebesgue space on $(X,\cal{B},m)$, where $\cal{B}$ denotes the Borel measurable subsets of $G$. Let $C_0(G)$ denotes the space of continuous functions on $G$ that vanish at infinity. Let $M(G)$ denote the space of finite, regular, Borel measures (Radon measures) on $G$. If $\chi \in \hat{G}$ and $f\in L^1(G)$, then the Fourier tranform of $f$ is defined by,  
\[
\hat{f}(\chi)=\int_G f(g)\cl{\chi(g)}\,dm(g).
\]

We will say that the dual of $\hat{G}$ is ordered if there exists a semigroup $P\subseteq \hat{G}$ such that $\hat{G}$ is the  disjoint union of $P$, $P^{-1}$ and $e$. We will denote by $H^p(G)$ the space of all elements in $L^p(G)$ such that $\widehat{f}|_{P^{-1}}=0$. 

Define $\ast:L^p(G)\to L^p(G)$ for $1\leq p\leq \infty$ by $f^\ast(g)=\cl{f(g^{-1})}$. The map $\ast$ is additive ($(f+g)^\ast = f^\ast+g^\ast$), conjugate-linear ($(\lambda f)^\ast = \cl{\lambda}f^\ast$), involutive ($f^{\ast\ast} = f$) and, since the group is unimodular, isometric ($\norm{f}=\norm{f^\ast}$). 

A simple calculation shows us that, 
\begin{align*}
\widehat{f^\ast}(\chi)&=\int_G f^\ast(g)\cl{\chi(g)}\,d\mu(g)=\int_G \cl{f(g^{-1})\chi(g)}\,d\mu(g)\\
&=\cl{\int_G f(g^{-1})\chi(g)}\,d\mu(g)=\cl{\int_G f(g)\chi(g^{-1})}\,d\mu(g)\\
&=\cl{\widehat{f}(\chi)}.
\end{align*}
It follows that $f=f^\ast$ if and only if the transform of $f$ is real valued and that the spaces $H^p$ are mapped into themselves by $\ast$. We will denote by $L^p_r(G)$ the set of elements of $L^p(G)$ whose Fourier transforms are real-valued. We will define $H^p_r(G):=L^p_r(G)\cap H^p(G)$ which is the set of elements in $H^p(G)$ whose Fourier transform is real valued.The map $\ast$ also is well-defined on $C_0(G)$ and is isometric.

Let $\Phi:L^p(G)\to L^p(G)$ be defined by,
\[\Phi(f)=\frac{f+f^\ast}{2}.\] 
The map $\Phi$ is a contractive projection from $L^p(G)$ onto $L^p_r(G)$. When restricted to $H^p(G)$, $\Phi$ carries $H^p(G)$ onto $H^p_r(G)$. The map $\Phi$ fixes any element in $L^p_r(G)$ and consequently $\Phi^2=\Phi$. If $L^p(G)$ is viewed as a real Banach space, then $L^p_r(G)$ is a closed subspace and $\Phi$ is a contractive linear projection. Every element $f\in L^p(G)$ can be written uniquely as $f=g+ih$ where $g,h\in L^p_r(G)$, this can be seen from  
\[
f=\frac{f+f^\ast}{2}+i\frac{f-f^\ast}{2i}=\Phi(f)+(f-\Phi(f)).
\]
Similar statements are true for $C_0(G)$. The most important property of $\Phi$ for our work is the fact that it is a module map. In particular if $f\in L^p_r(G)$ and $g\in L^q(G)$, where $p^{-1}+q^{-1}=1$, then $\Phi(fg)=f\Phi(g)$. This is easily seen from the definition.

For measures $\mu\in M(G)$ we define $\mu^\ast (E)=\cl{\mu(E^{-1})}$ and note that this defines a isometric map of $M(G)$ onto itself. As a functional on $C_0(G)$, 
\[
\mu^\ast(f)=\cl{\nu(f^\ast)}=\cl{\int_G f^\ast(g) d\mu(g) }.
\]
$\mu$  has a Fourier tranform defined by,
\[
\mu(\chi)=\int_G \chi(g)\,d\mu,
\]
and this agrees with the tranform on $L^1(G)$ under the natural identification of $L^1(G)$ as a subspace of $M(G)$.
Here too $\Phi$ is defined and has the properties mentioned above. 

\begin{lemma}\label{integrate}
If $f\in L^p(G)$ and $g\in L^q(G)$, with $p,q$ conjugate indices, then,
\[
\int f^\ast g=\cl{\int f g^\ast}.
\]
If $f\in C_0(G)$ and $\mu\in M(G)$, then,  
\[
\int_G f\,d\mu^\ast=\cl{\int_G f^\ast \,d\mu},
\]
for $f\in C_0(G)$ and $\mu\in M(G)$.
\end{lemma}

\begin{proof}
We have,
\begin{align*}
\int_G f^\ast(t)g(t)\,dt &= \int_G \cl{f(t^{-1})}g(t)\,dt =\cl{\int_G f(t^{-1})\cl{g(t)}\,dt }\\
&=\cl{\int_G f(t)\cl{g(t^{-1})}\,dt }=\cl{\int_G f(t)g^\ast(t)\,dt}.
\end{align*}
Note that the unimodularity of the group was used in the change of variables. The second statement about measures was mentioned earlier. By the definition of $\mu^\ast$ the result is true for simple functions and follows easily for $C_0(G)$ by an approximation argument.
\end{proof}

For $L^\infty(G)$ and $M(G)$ there is an additional property of $\Phi$ that is important.

\begin{prop}
Let $\Phi:L^\infty(G)\to L^\infty(G)$ be defined by $\Phi(f)=\frac{f+f^\ast}{2}$ and define it similarly for $M(G)$. In both cases $\Phi$ is $\wk$ continuous.
\end{prop}

\begin{proof}
We need only to check that $f\mapsto f^\ast$ is $\wk$ continuous. This follows from Lemma \ref{integrate}.
\end{proof}

We will now describe complexifcations for $L^p(G)$. These definitions apply equally well to the spaces $H^p(G)$, $M(G)$ or $C_0(G)$. If $\cal{M}$ is a subspace on $L^p_r(G)$ we define the complexification of $\cal{M}$ to be the subset of $L^p(G)$ given by 
\[\cal{M}_\bb{C}:=\cal{M}+i\cal{M}=\{f+ig\,:\,f,g\in \cal{M}\}.\]

\begin{lemma}
If $\cal{M}$ be a closed subspace of $L^p_r(G)$, $C_{0,r}(G)$ or $M_r(G)$, then $\cal{M}_\bb{C}$ is closed in $L^p(G)$, $C_0(G)$ or $M(G)$ respectively. If $p=\infty$ and $\cal{M}$ is a subalgebra of $L^\infty_r(G)$ or $C_{0,r}(G)$, then $\cal{M}_\bb{C}$ is an algebra. If $p=\infty$ and $\cal{M}$ is $\wk$ closed in $L^\infty_r(G)$, then $\cal{M}_\bb{C}$ is $\wk$ closed in $L^\infty$. If $\cal{M}$ is a $\wk$ closed subspace of $M_r(G)$, then $\cal{M}_\bb{C}$ is $\wk$-closed in $M(G)$.  
\end{lemma} 

\begin{proof}
We will prove this for $L^p(G)$, the other cases are similar. The fact that it is a subspace or an algebra is straightforward. To see that it is closed let $f_n+ig_n\to h$. Since $\Phi$ is continuous we have that $\Phi(f_n+ig_n)=f_n\to \Phi(h)$ and so $f_n+ig_n\to \Phi(h)+(h-\Phi(h))\in \cal{M}+i\cal{M}$. The statements about the $\wk$-closed subspaces follow from the $\wk$ continuity of $\Phi$.
\end{proof}

For $\chi\in L^\infty_r(G)$ we have,
\begin{equation}
\label{realtrans}\widehat{f\chi}(\eta)=\int_G f(t)\chi(t)\eta(t^{-1})\,dt=\int_G f(t)(\chi^{-1}\eta)(t^{-1})\,dt=\hat{f}(\chi^{-1}\eta).
\end{equation}
If $\int_G f$ and $\int_G g$ are both real then 
\begin{equation}
\label{realint}\int_G fg=\frac{1}{4}\left[\int(f+g)^2-\int_G (f-g)^2\right],
\end{equation} 
is real. If $p^{-1}+q^{-1}=1$, $f\in L^p_r(G)$, $g\in L^q_R(G)$  and  $\chi\in\hat{G}$, then an application of \eqref{realtrans} shows that $g\chi\in L^q_r(G)$. The equation \eqref{realint} now yields  $\int_G fg\chi\in \bb{R}$ and this proves that $fg$ is in $L^1_r(G)$. 

\begin{lemma}
Let $1\leq p<\infty$. The dual of $L^p_r(G)$ is isometrically isomorphic to $L^q_r(G)$. 
\end{lemma}

\begin{proof}
It is well-known that if $X$ is a real Banach space and $X_\bb{C}$ is its complexification, then the dual of $X$ is isometrically isomorphic to the space $\{\Re (\phi)\,:\, \phi\in (X_\bb{C})^\ast\}$. Let $f\in L^q(G)$ and let $\phi_f$ denote the functional it induces on $L^1(G)$. Let $g\in L^p_r(G)$ and consider,
\begin{align*}
2(\Re (\phi_f) )(g)&= \phi_f(g)+\cl{\phi_f}(g)=\int_G fg+\int_G \cl{fg}\\
&=\int_G fg+\int\cl{fg}=\int_G fg +\int_G \cl{fg^\ast}\\
&=\int_G fg +\int_G f^\ast g=2\int_G \Phi(f)g
\end{align*}
Thus, $L^p_r(G)^\ast=\{\Re (\phi_f) \,:\, f\in L^q(G)\}=\Phi(L^q(G))=L^q_r(G)$.
\end{proof}

\begin{lemma}
The dual of $C_{0,r}(G)$ is isometrically isomorphic to $M_r(G)$.
\end{lemma}

\begin{proof}
As in the previous theorem the dual of $C_{0,r}(G)$ is identified with the the real parts of functionals on $C_0(G)$. As before if $\mu\in M(G)$, then, 
\begin{align*}
2(\Re \mu )(f)&= \mu(f)+\cl{\mu(f)}=\int_G f \,d\mu+\cl{\int_G f\, d\mu}\\
&=\int_G f d\,\mu + \cl{\int_G f^\ast \,d\mu} =\int_G f \,d\mu +\int_G f \,d\mu^\ast\\
&=2\int_G f\,d\Phi(\mu).
\end{align*}
Therefore, $C_{0,r}(G)^\ast=M_r(G)$.
\end{proof}

Similar duality results hold for the real Hardy spaces and we mention one of these. It is well-known that the dual of $H^1$ can be identified with the space BMOA. Every function in BMOA has a Fourier expansion and we see that the dual of $H^1_r(\bb{T})$ can be identified with $\mathrm{BMOA}_r$. 

For the remainder of the paper we will consider $L^\infty(G)$ endowed with its $\wk$ topology. Thus all statements regarding closures and convergence in $L^\infty(G)$ will be with respect to the $\wk$ topology. The same assumption will be made about $H^\infty(G)$ and $M(G)$.

For us the most interesting cases will be $G=\bb{T}$, $G=\bb{R}$ and in these cases $\hat{G}=\bb{Z}$ and $\hat{G}=\bb{R}$, respectively. We note that in both cases the dual group is ordered. 

\section{Invariant Subspaces and Maximal Algebras}

 We will denote by $z$ the identity function on $\bb{T}$. Let $1\leq p\leq \infty$. A function $\phi$ in $H^p(\bb{T})$ is called outer if the set $\{\phi z^n\,:\,n\geq 0\}$ spans a dense subspace of $H^p(\bb{T})$.

\begin{lemma}
Let $\phi$ be an outer function in $H^p(\bb{T})$, then $\phi^\ast$ is outer.
\end{lemma}

\begin{proof}
Let $p$ be finite and let $q$ be the index conjugate to $p$. Assume that $f\in L^q(\bb{T})$ and that $f$ is in the annihilator of the subspace spanned by $z^n \phi^\ast$, $n\geq 0$. By the above lemma we have that $f^\ast$ is in the annihilator of $(z^n \phi^\ast )^\ast=z^n\phi$. Since $\phi$ is outer this imples $f^\ast=0$ and so $f=0$. Therefore, $\phi^\ast$ is outer. A similar proof works for $p=\infty$ by considering $f$ in the preannihilator.
\end{proof}

As far as inner functions are concerned we note that $\mod{\phi}=1$ a.e. on $\bb{T}$ if and only if $\mod{\phi^\ast}=1$ a.e. on $\bb{T}$. Hence if $\phi$ is unimodular, then $\phi^\ast$ is unimodular. 

\begin{thm}
Let $f\in H^1_r(\bb{T})$, then there exists $\phi$ inner and $u$ outer such that $\phi=\phi^\ast$, $u=u^\ast$ and  $f=\phi u$. The factorization is unique up to sign change.
\end{thm}

\begin{proof}
Factor $f$ as $f=\phi u$ with $\phi$ inner and $u$ outer. Note that $\phi^\ast u^\ast =f^\ast = f=\phi u$. Multiply $\phi$ by a unimodular scalar so that it has a non-zero real Fourier coefficient. Now $\phi=\lambda \phi^\ast$ and $u=\cl{\lambda} u^\ast$. Comparing coefficients we get $\lambda=1$. We know that the factorization is unique up to unimodular scalar factors. Imposing the extra condition that $\phi=\phi^\ast$ implies that the only permissible scalars are $\pm 1$.
\end{proof}

Our next result is the analogue of Riesz factorization.

\begin{thm}
If $f\in H^1_r(\bb{T})$, then there exists $f_1,f_2\in H^2_r(\bb{T})$ such that $f=f_1f_2$ and $\norm{f}_1=\norm{f_1}_2^2=\norm{f_2}_2^2$. 
\end{thm}

\begin{proof}
We can factor $f=\phi u$  with $\phi$ inner and $u$ outer and $\phi,u\in H^1_r(\bb{T})$. Set $g=u^{1/2}$ which clearly has the property that $g^\ast=g$. Alternatively comparing the coefficients of $g^2$ and  $u$ yields the fact that $g$ is in $H^2_r(\bb{T})$. Setting $f_1=\phi g$ and $f_2=g$ yields the factorization.
\end{proof}

We will describe the closed shift invariant subspaces of $L^p_r$ for $1\leq p<\infty$, the $\wk$ closed shift invariant subspaces of $L^\infty_r(\bb{T})$ and the closed shift invariant subspaces of $C_r(\bb{T})$.

Let $E\subset \bb{C}$. We will denote by $E^\ast$ the set of complex conjugates of points in $E$. By shift invariant  we mean a subspace $\cal{M}$ of $L^p(\bb{T})$ or $L^p_r(\bb{T})$ such that $z\cal{M}\subseteq \cal{M}$. We call a shift invariant subspace \textit{simply invariant} if $z\cal{M}\not=\cal{M}$ and \textit{doubly invariant} if $z\cal{M}=\cal{M}$. 

\begin{thm}
Let $\cal{M}$ be a closed shift invariant subspace of $L^p_r(\bb{T})$. The subspace $\cal{M}$ is simply invariant if and only if $\cal{M}=\phi H^p_r(\bb{T})$ for some unimodular function $\phi$ with $\phi=\phi^\ast$. The subspace $\cal{M}$ is doubly invariant if and only if there exists a measurable subset $E$ of the circle $\bb{T}$, such that $E=E^\ast$, and  $\cal{M}=\chi_EL^p_r(\bb{T})$.
\end{thm}

\begin{proof}
Begin by noting that $\cal{M}_\bb{C}\subseteq L^1(\bb{T})$ is invariant under $z$ and so is either of the form $\phi H^p$ for some unimodular $\phi$ or of the form $\chi_E L^1(\bb{T})$ for some measurable subset $E\subseteq \bb{T}$. If $\cal{M}_\bb{C}=\phi H^p$ we may multiply $\phi$ by a unimodular scalar so that it has at least one real non-zero Fourier coefficient.

Note that $\cal{M}^\ast:=\{f^\ast\,:\,f\in \cal{M}\}=\cal{M}$ and consequently $\cal{M}_\bb{C}^\ast=\cal{M}_\bb{C}$. If $\cal{M}_\bb{C}=\phi H^1(\bb{T})$, then we have $\phi^\ast H^1(\bb{T})=\cal{M}_\bb{C}^\ast =\cal{M}_\bb{C}=\phi H^1(\bb{T})$ and so $\phi^\ast=\lambda \phi$ for some $\lambda\in \bb{T}$. Since $\phi$ was chosen to have at least one non-zero real Fourier coefficient we see that $\lambda = 1$ and $\phi=\phi^\ast$. If $\cal{M}_\bb{C}=\chi_E L^1(\bb{T})$, then similar reasoning shows that $\chi_E=\chi_E^\ast$ and it is easy to check this is equivalent to $E=E^\ast$.

To finish the proof observe that if $f=f^\ast$, then $\Phi(fg)=f\Phi(g)$ and so $\cal{M}=\Phi(\cal{M}_\bb{C})=\phi\Phi(H^1(\bb{T}))=\phi H^1_r(\bb{T})$ or  $\cal{M}=\Phi(\cal{M}_\bb{C})=\Phi(\chi_E L^1(\bb{T}))=\chi_EL^1_r(\bb{T})$.

\end{proof}

Let $\cal{S}$ be a subset of $C(\bb{T})$. We will denote by $Z(\cal{S})$ the set of common zeros on $\bb{T}$ of functions in $\cal{S}$. We will denote by $I(\cal{S})$ the ideal of functions that vanish on $Z(\cal{S})$. We begin with the following observation.

\begin{lemma}
If $\cal{S}\subseteq C(\bb{T})$, then $Z(\cal{S}^\ast)=Z(\cal{S})^\ast$ and $I(\cal{S}^\ast)=I(\cal{S})^\ast$.
\end{lemma}

\begin{proof}
If $\lambda\in Z(\cal{S}^\ast)$, then $f^\ast(\lambda)=0$ for all $f\in \cal{S}$ and hence $f(\cl{\lambda})=0$ for all $f\in \cal{S}$ and $\lambda \in Z(\cal{S})^\ast$. Therefore, $Z(\cal{S}^\ast)\subseteq Z(\cal{S})^\ast$. Now, $Z(\cal{S})^\ast =Z(\cal{S}^{\ast\ast})^\ast \subseteq Z(\cal{S}^\ast)$ and we are done. The second claim follows easily from the fact that $f^\ast|E=0$ if and only if $f|E^\ast=0$. 
\end{proof}

\begin{thm}
Let $\cal{M}$ be a closed subspace of $C_r(\bb{T})$ which is invariant under multiplication by $z$. If $z\cal{M} \not= \cal{M}$, then there exists a unimodular function $\phi$ such that $\phi=\phi^\ast$ and $\cal{M}=\phi H^\infty_r(\bb{T})\cap I(\cal{M})$. If $z \cal{M}=\cal{M}$, then $\cal{M}=I(\cal{M})$. 
\end{thm}

\begin{proof}
By the main result of \cite{HS} the closed subspace $\cal{M}_\bb{C}$ has the form $\phi H^\infty(\bb{T})\cap I(\cal{M}_\bb{C})$ or is equal to $I(\cal{M}_\bb{C})$. Since $\cal{M}_\bb{C}^\ast =\cal{M}_\bb{C}$, we have that $I(\cal{M}_\bb{C})^\ast =I(\cal{M}_\bb{C})=I(\cal{M})$. In the case where $\cal{M}_\bb{C}=\phi H^\infty(\bb{T})\cap I(\cal{M}_\bb{C})$, the uniqueness of $\phi$ up to multiplication by a unimodular scalar implies as before that we may choose $\phi=\phi^\ast$. Applying $\Phi$ we get,
\[\cal{M}=\Phi(\phi H^\infty(\bb{T})\cap I(\cal{M}_\bb{C}))=\phi\Phi(H^\infty(\bb{T}))\cap I(\cal{M}) = \phi H^\infty_r(\bb{T})\cap I(\cal{M}).\]
\end{proof}

If $K$ is a compact subset of the circle $\bb{T}$, let $M(K)$ denote the set of measures in $M(\bb{T})$ that are supported on $K$.

\begin{thm}
Let $\cal{M}$ be a $\wk$ closed invariant subspace of $M_r(\bb{T})$. Then there exists a unimodular $\phi\in L^\infty_r(\bb{T})$ and $K$ of Lebesgue measure $0$ such that $K=K^\ast$ and $\cal{M}=\phi H_{0,r}^1(\bb{T})+M(K)$.
\end{thm}

As in the classical case, this can be deduced from the corresponding result about $C_r(\bb{T})$. We will deduce it directly from the classical result.

\begin{proof}
The results of \cite{HS} show that  $\cal{M}_\bb{C}$ is of the form $\phi H_0^1(\bb{T})+M(K)$ and $K$ is a closed subset of $\bb{T}$ of Lebesgue measure 0. The set $K$ is uniquely determined by $\cal{M}_\bb{C}$ and $\phi$ is unique upto multipication by unimodular scalars. By applying $\ast$ to $\cal{M}_\bb{C}$ and noting that it is unchanged we get that $\phi^\ast H_0^1(\bb{T})+M(K)=\phi H_0^1 (\bb{T})+M(K^\ast)$. It follows that $K=K^\ast$ and that we may choose $\phi$ such that $\phi=\phi^\ast$. Applying the projection $\Phi$ we obtain $\cal{M}=\phi H_{0,r}^1(\bb{T})+M(K)$.
\end{proof}

By a duality argument we can also show that the $\wk$ closed invariant subspaces of $\mathrm{BMOA}_r$ have the same form as in the classical case. 

Our next application is Wermer's maximality theorem and the structure of closed ideals in $A_r(\bb{D})$.

\begin{thm}
If $\cal{B}$ is a closed algebra such that  $A_r(\bb{D})\subseteq \cal{B} \subseteq C_r(\bb{T})$, then $\cal{B}=A_r(\bb{D})$ or $\cal{B}=C_r(\bb{T})$.
\end{thm}

\begin{proof}
Let $\cal{B}_\bb{C}=\cal{B}+i\cal{B}$ and note that $\cal{B}_\bb{C}$ is a uniformly closed subalgebra of $C(\bb{T})$ containing $A(\bb{D})$. Hence $\cal{B}_\bb{C}=A(\bb{D})$ or $B_\bb{C}=C(\bb{T})$. Applying $\Phi$ we see that $\cal{B}=\Phi(B_\bb{C})=A_r(\bb{D})$ or $\cal{B}=C_r(\bb{T})$.
\end{proof}

An identical argument yields the following. The proof in the classical case may be found in \cite{KH}.

\begin{thm}
If $\cal{B}$ is a $\wk$ closed subalgebra of $L_r^\infty(\bb{T})$ containing $H^\infty_r(\bb{T})$, then $\cal{B}=H^\infty_r(\bb{T})$ or $\cal{B}=L^\infty_r(\bb{T})$. 
\end{thm}

\begin{thm}
If $\cal{J}$ is a closed ideal in  $A_r(\bb{D})$, then $\cal{J}=\phi H^\infty_r(\bb{T}) \cap I(\cal{J})$ for some inner function $\phi$ in $H^\infty_r(\bb{T})$. If $\cal{J}$ is maximal, then there exists a $\lambda\in \cl{\bb{D}}$ such that $\cal{J}=\{f\in A_r\,:\,f(\lambda)=f(\cl{\lambda})=0\}$.
\end{thm}

\begin{proof}
A result of Rudin, see \cite{KH}, tells us that the closed ideal $\cal{J}_\bb{C}$ is of the form $\phi H^\infty(\bb{T}) \cap I(\cal{J}_\bb{C})$. As in the case of invariant subspaces of $C_r(\bb{T})$ it follows that $\cal{J}=\phi H^\infty_r(\bb{T})\cap I(\cal{J})$. It is easy now to see that the maximal ideals are of the form given above.  
\end{proof}

For the real line the earliest and most widely known invariant subspace results are due to Wiener and Lax. Wiener's theorem classifies the closed subspaces of $L^2(\bb{R})$ that are invariant for multiplication by $e^{ixt}$ for all $x\in \bb{R}$. These can also be viewed, by applying the Fourier transform, as the translation invariant subspaces of $L^2(\bb{R})$. The theorem of Lax gives the structure of the subspaces invariant for multiplication by $e^{ixt}$ for $x>0$. The following theorem is the corresponding ``real'' result.

\begin{thm}[Wiener \& Lax]
Let $\cal{M}$ be a closed subspace of $L^2_r(\bb{R})$. If $\cal{M}$ is invariant under multiplication by $e^{ixt}$ for all $x\in \bb{R}$, then there exists a measurable set $E$ such that $E=E^\ast$,
\[\{\hat{f}\,:\,f\in \cal{M}\}=\chi_EL^2_r(\bb{R}).\]
If $\cal{M}$ is invariant under multiplication by $e^{ixt}$ for all $x>0$, but not for all $x\in \bb{R}$, then there exists a unimodular function $\phi$ defined on $\bb{R}$ such that $\cal{M}=\phi H^2_r(\bb{R})$ and $\phi^\ast=\phi$.
\end{thm}

Here $\phi$ is unimodular and consequently does not belong to any of the $L^p$ spaces. However the definition of $\ast$ still makes sense and we have $\phi(t)=\cl{\phi(-t)}$ for all $t\in \bb{R}$. The set $E$ described above  now symmetric about the origin.

We remark that these projection techniques continue to hold in the setting of compact groups with ordered duals as defined by Helson in \cite{HH2}. The idea of using an expectation to classify invariant subspaces seems to have appeared first in a paper of Abrahamse and Douglas \cite{AD}.

\section{F. \& M. Riesz Theorem and Szeg\"{o}'s Theorem}\label{measures}

The classical result of F. \& M. Riesz shows that an analytic measure on the circle is automatically absolutely continuous with respect to Lebesgue measure. The next result is the analogue of this result for the real Hardy spaces.

\begin{thm}
Let $\mu\in M_r(\bb{T})$ with $\hat{\mu}(n)=0$ for all $n<0$. There exists a function $f\in H^1_r(\bb{T})$ such that $d\mu =f d\theta$ and $\norm{\mu}=\norm{f}_{1}$.
\end{thm}

\begin{proof}
We know that there exists an $f\in H^1(\bb{T})$ such that $d\mu =f d\theta$ with $\norm{\mu}=\norm{f}_1$. It follows that $\hat{\mu}(n)=\hat{f}(n)\in \bb{R}$ and so $f\in H^1_r(\bb{T})$.
\end{proof}

On the space of regular Borel measures $M(\bb{T})$ the action of $\ast$ is given by $\mu^\ast(E)=\cl{\mu(E^\ast)}$. It is clear from this that if $\mu$ is singular, then so is $\mu^\ast$ and if $\mu$ is absolutely continuous with respect to Lebesgue measure and $d\mu=f d\theta$, then $d\mu^\ast =f^\ast d\theta$.

\begin{thm}
Let $\mu\in M(\bb{T})$ and denote by $\mu_a$ and $\mu_s$ the absolutely continuous and singular parts of $\mu$. If $\mu\in M_r(\bb{T})$, then $\mu_a,\mu_s\in M_r(\bb{T})$. 
\end{thm}

\begin{proof}
Note that $\mu_a+\mu_s=\mu=\mu^\ast=\mu_a^\ast+\mu_s^\ast$. Since the absolutely continuous and singular parts are uniquely determined we see that $\mu_a=\mu_a^\ast$ and $\mu_s=\mu_s^\ast$. 
\end{proof}

\begin{thm}
Let $\mu\in M_r(\bb{T})$, let $\mu_a$ be its absolutely continuous part and $w=d\mu_a/d\theta \in L^1_r$. We have, 
\begin{equation}
\label{szego1}\inf \left\{\int_{\bb{T}}\mod{1-f}^2\,d\mu\,:\,f\in A_{r,0}(\bb{D})\right\}=\exp\int \log w\, d\theta.
\end{equation}
\end{thm}

\begin{proof}
The integral on the right of \eqref{szego1} is equal to 
\[\inf\left\{\int_\bb{T}\mod{1-f}^2w\,d\theta \,:\,f\in A_0(\bb{D})\right\},\]
by Szeg\"{o}'s theorem. Consider the map $f\mapsto f^\ast$ on $L^2(d\mu_a)$. We have,
\[\int_{\bb{T}} \mod{f^\ast}^2w\,d\theta=\int_{\bb{T}}(\mod{f}^2)^\ast w\,d\theta=\cl{\int_\bb{T} \mod{f}^2w^\ast\,d\theta}=\int_\bb{T}\mod{f}^2w\,d\theta.\]
The last step used the fact that $w=w^\ast$ and the fact that $w$ is non-negative. Hence the map $\Phi(f)=\frac{f+f^\ast}{2}$ is a contractive projection on $L^2(d\mu_a)$ and consequently,
\[\norm{1-f}_{L^2(d\mu_a)}^2\leq \norm{1-\Phi(f)}^2_{L^2(d\mu_a)}.\]
Using this last inequality and the fact that $\Phi$ projects $A_0(\bb{D})$ onto $A_{r,0}(\bb{D})$ we get,
\begin{align*}
& \inf \{\norm{1-f}_{L^2(d\mu_a)}^2\,:\, f\in A_0(\bb{D})\}\\
\leq & \inf \{\norm{1-f}_{L^2(d\mu_a)}^2\,:\, f\in A_{r,0}(\bb{D})\}\\
=&\inf \{\norm{1-\Phi(f)}_{L^2(d\mu_a)}^2\,:\, f\in A_0(\bb{D})\}\\
\leq &\inf \{\norm{1-f}_{L^2(d\mu_a)}^2\,:\, f\in A_0(\bb{D})\},
\end{align*}
which establishes the equality in \eqref{szego1}.
\end{proof}

\section{Interpolation results}\label{misc}

In this section we prove a Nevanlinna-Pick interpolation theorem for the real Hardy algebra $H^\infty_r$. We also look at the Carath\'eodory-Fej\'er interpolation problem.

Given $n$ points $z_1,\ldots,z_n\in \bb{D}$ and $n$ complex scalars $w_1,\ldots,w_n\in\bb{C}$ the classical Nevanlinna-Pick theorem provides a necessary and sufficient condition for the existence of a holomorphic map $f:\bb{D}\to \bb{D}$ such that $f(z_j) = w_j$ for $j=1,\ldots,n$. The condition is the positivity of the $n\times n$ Pick matrix 
\[\left[\frac{1-w_i\cl{w_j}}{1-z_i\cl{z_j}}\right]_{i,j=1}^n.\]
We now give a condition that guarantees the existence of a function $f\in H^\infty_r$, such that $\norm{f}_\infty\leq 1$ and $f(z_j) = w_j$ for $j=1,\ldots,n$. We denote by $r$ the number of points in the set $\{z_1,\ldots,z_n\}$ which lie on the real axis. We also reorder the points so that $z_1,\ldots,z_r\in\bb{R}$. Suppose that $f\in H^\infty_r$ and $f(z_j) = w_j$ for $j=1,\ldots,n$. It follows that 
\[w_j = f(z_j) = f^\ast(z_j) = \cl{f(\cl{z_j})}.\]
Hence, $w_j \in \bb{R}$ for $j=1,\ldots,r$ and $f(\cl{z_j}) = \cl{w_j}$ for $j=r+1,\ldots,n$. Conversely assume that $w_1,\ldots,w_r\in\bb{R}$ and suppose there exists a function $f\in H^\infty$ such that $f(z_j) = w_j$ for $j=1,\ldots,n$ and $f(\cl{z_j})=\cl{w_j}$ for $j=r+1,\ldots,n$. The function $g=\Phi(f)\in H^\infty_r$, $\norm{g}_\infty\leq \norm{f}_\infty$ and 
\[g(z_j) = \dfrac{f(z_j)+\cl{f(\cl{z_j})}}{2} = w_j.\]
Let $s$ be the number of points for which both $z_j$ and $w_j$ are real. We know that $0\leq s\leq r$ and suppose that the Pick matrix for the points $z_1,\ldots,z_n,z_{s+1},\ldots,z_r,\cl{z_{r+1}},\ldots,\cl{z_n}\in\bb{D}$, $w_1,\ldots,w_n,\cl{w_{s+1}},\ldots,\cl{w_n}\in\bb{D}$ is positive. In particular, the $2\times 2$ submatrix corresponding to $w_j$, $\cl{w_{j}}$ for $s<j\leq r$ is a Pick matrix $P$ for points of the form $z_j,z_j\in\bb{D}$, $w_j,\cl{w_j}\in\bb{D}$. The positivity of this matrix gives
\[
\begin{bmatrix}
\frac{1-\mod{w_j}^2}{1-\mod{z_j}^2} & \frac{1-w_j^2}{1-\mod{z_j}^2} \\ \frac{1-\cl{w_j}^2}{1-\mod{z_j}^2} & \frac{1-\mod{w_j}^2}{1-\mod{z_j}^2}
\end{bmatrix}\geq 0.\]
Taking the determinant and simplifying we get,
\[\Re (w_j^2)\geq \mod{w_j}^2, \]
which shows that $w_j\in\bb{R}$ for $s<j\leq r$. If $P$ is positive, then $w_j = \cl{w_j}$ for $j=1,\ldots,r$ and  Nevanlinna-Pick theorem implies there exists $f\in H^\infty$ such that $\norm{f}_\infty\leq 1$ with $f(z_j) = w_j$ for $j=1,\ldots,n$ and $f(\cl{z_j})=\cl{w_j}$ for $j=r+1,\ldots,n$. It follows that $g=\Phi(f)$ solves the interpolation problem in $H^\infty_r$. These arguments prove the following result. 

\begin{thm}[Nevanlinna-Pick interpolation for $H^\infty_r$.]
Let $z_1,\ldots,z_n\in \bb{D}$ and let $w_1,\ldots,w_n\in \bb{C}$. Let $r,s$ and $P$ be as described in the above paragraph. There exists a function $f\in H^\infty_r$ such that $f(z_j) = w_j$ if and only if $P$ is positive.
\end{thm}

It would seem natural to consider the Carath\'eodory-Fej\'er interpolation problem as well. In the Carath\'eodory-Fej\'er problem the first $n+1$ Taylor coefficients $a_0,\ldots,a_n$ are specified. A necessary and sufficient condition for the existence of a function $f\in H^\infty$ such that $\norm{f}_\infty\leq 1$ and $f(z) = a_0+\ldots+a_nz^n+z^{n+1}g(z)$ for some $g\in H^\infty$ is that the T\"oplitz matrix 
\[
\begin{bmatrix}
a_0 &  & 0\\ 
\vdots & \ddots &  \\
a_n & \cdots  & a_0
\end{bmatrix}
\]
be a contraction. If we start with real scalars $a_0,\ldots,a_n$ and apply the interpolation result we see that there exists $f,g\in H^\infty$ such that $\norm{f}_\infty\leq 1$ and $f(z)=a_0+\ldots+a_nz^n+z^{n+1}g(z)$. It is easy to check that $\Phi(f)(z) =  a_0+\ldots+a_nz^n+z^{n+1}\Phi(g)$ and so $\Phi(f)$ solves the Carath\'eodory-Fej\'er interpolation problem in $H^\infty_r$.

\section{Extensions of Real Function Theory}

In general one would like to study the analogous questions about invariant subspaces for algebras besides $H^\infty$. General uniform algebras provide one possible extension of the theory. In studying general uniform algebras the correct notion of real is no longer clear. Another possble extension is to consider subalgebras of $H^\infty$. We show how the theory extends to one specific case that has received some attention recently \cite{DPRS}. 

The process of taking real parts of non-constant analytic functions does not keep us in the space of analytic functions. This makes the idea of looking at the functions with real coefficients interresting since it provides a naturally occuring example of a real vector space of complex analytic functions. We will also look at the codimension one subalgebra of $H^\infty_r$, generated by $z^2,z^3$ and show that the projection argument allows to extend the invariant subspace results in \cite{DPRS}. In \cite{DPRS} a classification is given for subspaces of $H^2$ that are invariant for multiplication by $z^2$ and $z^3$ but not for $z$. Such a subspace can be written in a unique way as $\phi([1+\alpha z]\oplus z^2H^2)$. where $\phi$ is an inner function and $[1+\alpha z]$ denotes the span in $H^2$ of the vector $1+\alpha z$. Let $\psi$ be an inner function and $\beta\in \bb{C}$ with $\phi([1+\alpha z]\oplus z^2H^2)=\psi([1+\beta z]\oplus z^2H^2)$. We may assume that $\phi(0)\not = 0$ and write $\psi z^2 = \phi(\lambda(1+\beta )+z^2h)$ for some $\lambda\in \bb{C}$ and $h\in H^2$. Evaluating at $0$ we get $\phi(0)\lambda = 0$ and so $\lambda = 0$ and $\psi = \phi h$. If $\psi(0)=0$, then $h(0)=0$. If we write $(1+\alpha z)= \psi (\mu(1+\beta z)+z^2 q)$, then we see that $z$ divides the function $\psi (\mu(1+\beta z)+z^2 q)$ but not $(1+\alpha z)$. Hence, $\psi(0)\not = 0$. It follows now that there exists $g$ such that  $\phi = \psi g$ and so $\phi$ and $\psi $ are uniquely determined upto a unimodular factor. The fact that $\alpha = \beta$ follows easily. 

Having established this we can now prove the an analogue of the result in \cite{DPRS} for the functions $f\in H^\infty$ such that $f'(0)=0$ and $f$ has real Fourier coefficients.

\begin{prop}
Let $\cal{M}$ be a closed subspace of $H^2_r$ such that $z^2\cal{M},z^3\cal{M}\subseteq \cal{M}$ but $z\cal{M}\not\subseteq \cal{M}$. There exists a real scalar $c$ and an inner function $\phi \in H^\infty_r$ such that $\cal{M}=\phi([1+cz]\oplus z^2H^2_r)$.
\end{prop}

\begin{proof}
The result follows from the projection arguments once we note that $\cal{M}+i\cal{M}=\phi ([1+\alpha z]\oplus z^2H^2)$ and that $\phi$ can be chosen with real coefficients by the uniqueness result in the previous paragraph. Hence $\Phi(\phi([1+\alpha z]\oplus z^2H^2))=\phi([1+\Re(\alpha)]\oplus z^2H^2_r)$.
\end{proof}

Finally, we mention a minor open question. In \cite{MS} and \cite{JR} the invariant subspace result is established by a direct argument with no knowledge of Beurling' s theorem. Is it possible to deduce Beurling's theorem from the analogous real Fourier coefficients result?

The obvious approach does not seem to work. Let $\cal{M}$ be a closed invariant subspace of $H^2$. Assume in addition that $\cal{M}=\cal{M}^\ast$. For every $f\in \cal{M}$ we have that $f+f^\ast\in \cal{M}$. By modifying $f$ to have at least one real Fourier coefficient we see that $\cal{N}=\cal{M}\cap H^2_r$ is a non-trivial closed invariant subspace of $H^2_r$. From the real result we see that $\cal{N}=\phi H^2_r$ for some real inner function $\phi$. Our assumption that $\cal{M}=\cal{M}^\ast$ implies that for $f\in\cal{M}$, $f=\dfrac{f+f^\ast}{2}+i\dfrac{f-f^\ast}{2i}\in \cal{N}+i\cal{N}$. The inclusion $\cal{N}+i\cal{N}\subset \cal{M}$ is straightforward and so $\cal{M}=\phi H^2_r+i\phi H^2_r=\phi H^2$. Hence a shift invariant subspace $\cal{M}$ of $H^2$ with the property that is also invariant under $\ast$ has the form $\phi H^2$ for a real inner function $\phi$. 

For a general invariant subspace $\cal{M}$ we can always form the smallest $\ast$-invariant space that contains it, $\cal{M}+\cal{M}^\ast$. We see that $\cal{M}\subseteq \cal{M}+\cal{M}^\ast = \phi H^2$ for some real inner function $\phi$. However, the inclusion can be strict and the $\phi$ can be 1. Consider for example the space $\cal{M}_a$ of functions that vanish at $a$ where $a\not\in\{1,-1,i,-i\}$. The subspace $\cal{M}_a=\phi_aH^2$ where $\phi_a$ is the simple Blaskche factor for the point $a$. The inner function $\phi_a\not\in H^2_r$ and since $\cal{M}_a$ is of co-dimension 1 in $H^2$ we see that $\cal{M}_a+\cal{M}_a^\ast=H^2$.

\end{document}